\documentstyle[12pt,amsfonts,amscd]{amsart}
\newtheorem{Theorem}{Theorem}[section]
\newtheorem{Definition}[Theorem]{Definition}
\newtheorem{Proposition}[Theorem]{Proposition}

\newtheorem{Lemma}[Theorem]{Lemma}
\newtheorem{Corollary}[Theorem]{Corollary}
\theoremstyle{remark}
\newtheorem{Example}[Theorem]{Example}

\def\Aut{\operatorname{Aut}}
\def\dim{\operatorname{dim}}

\def\be{\begin{enumerate}}
\def\ee{\end{enumerate}}
\def\bT{\begin{Theorem}}
\def\eT{\end{Theorem}}
\def\bP{\begin{Proposition}}
\def\eP{\end{Proposition}}
\def\bD{\begin{Definition}}
\def\eD{\end{Definition}}
\def\bE{\begin{Example}}
\def\eE{\end{Example}}
\def\bL{\begin{Lemma}}
\def\eL{\end{Lemma}}
\def\bC{\begin{Corollary}}
\def\eC{\end{Corollary}}

\begin{document}
\title{Properties of Fixed Point Sets and a Characterization of the Ball in ${\Bbb C}^n$ }
\author{Buma L. Fridman and  Daowei Ma}
\begin{abstract} We study the fixed point sets of holomorphic self-maps of a bounded domain in ${\Bbb C}^n$. Specifically we investigate the least number of fixed points in general position in the domain that forces any automorphism (or endomorphism) to be the identity. We have discovered that in terms of this number one can give the necessary and sufficient condition for the domain to be biholomorphic to the unit ball. Other theorems and examples generalize and complete previous results in this area, especially the recent work of Jean-Pierre Vigu\'{e}.\end{abstract}
\address{ buma.fridman@@wichita.edu, Department of Mathematics,
Wichita State University, Wichita, KS 67260-0033, USA}
\address{ dma@@math.wichita.edu, Department of Mathematics,
Wichita State University, Wichita, KS 67260-0033, USA}
\maketitle \setcounter{section}{-1}
\section{Introduction}
Let $D\subset {\Bbb C}^n$ be a bounded domain. Below we consider two
families of self-maps of $D$. The first is the group $Aut(D)$ of
holomorphic automorphisms of $D$; the second is the set $H(D,D)$ of all
holomorphic maps from $D$ to $D$, i.e., the set of endomorphisms of $D$. 

\bD A set $K\subset D$ is called a determining subset of $D$ with
respect to $Aut(D)$ (or $H(D,D)$ resp.) if, whenever $g$ is automorphism
(resp. endomorphism) of $D$ such that $g(k)=k$ $\forall k\in K$, then $g$ is
the identity map of $D$.
\eD
The notion of a determining set was introduced earlier in a paper we wrote
with our collaborators Steven G. Krantz and Kang-Tae Kim [FK1]. In that
paper we attempted to find a higher dimensional analog of the following
result of classical function theory ([FF],[ES],[Mas],[PL],[Su]): if $f:M\rightarrow M$ is a conformal
self-mapping of a plane domain $M$ which fixes three distinct points then $%
f(\zeta )=\zeta $. 
\newline Determining sets have been further investigated in the following papers [FK2], [KK], [Vi1], [Vi2]. 

Let $W_s(D)$ denote the set of $s$-tuples $(x_1,\dots, x_s)$, where $x_j\in D$, such that $\{x_1,\dots, x_s\}$ is a determining set with respect to $Aut(D)$. Similarly, $\widehat{W}_s(D)$ denotes the set of $s$-tuples $(x_1,\dots, x_s)$ such that $\{x_1,\dots, x_s\}$ is a determining set with respect to $H(D,D)$.
So $\widehat{W}_s(D)\subseteq W_s(D)\subseteq D^s$.
We now introduce two numbers $s_0(D)$ and $\widehat{s}_0(D)$. In case $Aut(D)={id}$, $s_0(D)=0$, otherwise $s_0(D)$ is the least integer $s$, such that $W_s(D)\neq \emptyset $. The symbol $\widehat{s}_0(D)$ denotes the least integer $s$ such that $\widehat{W}_s(D)\neq \emptyset $. Hence, $s_0(D)\leq \widehat{s}_0(D)$.

In [FK1] we proved the inequality $s_0(D)\leq n+1$ for many (but not all)
bounded domains in ${\Bbb C}^n$. In [Vi1] J.-P. Vigu\'{e}, using a different
method proved this estimate for \textit{all }bounded domains in ${\Bbb C}^n$
. Furthermore in [Vi2] Vigu\'{e} proved the estimate $\widehat{s}%
_0(D)\leq n+1$ for all bounded domains in ${\Bbb C}^n$.

Both estimates are the best possible, since
for the unit ball $B^n\subset {\Bbb C}^n$, $s_0(B^n)=n+1$. In section 1 we
prove that the reverse is true: if $s_0(D)=n+1$  for a bounded domain $D\subset 
{\Bbb C}^n$, then $D$ is biholomorphic to the unit ball $B^n$. Obviously, $%
s_0(D)$ depends on how large $Aut(D)$ is: for a smaller group, we expect a
lower $s_0(D)$. This relationship is reflected below in Corollary \ref{in}.

If a positive integer $s\geq s_0(D)$, then $W_s(D)\neq \emptyset $, so there are $s$ points such that if an automorphism of 
$D$ fixes these points it will fix any point of $D$. Now the question arises
whether the choice of these $s$ points is generic. To make it more precise
we need to find out if $W_s(D)$ is open and everywhere dense in $D^s$. We consider this question in section 2. Refining and complementing the results of [FK1, Vi1, Vi2] we prove that $W_s(D)$ is open, and also dense if not empty. For the similar question related to $\widehat{W}_s(D)$ we provide examples to the contrary.

\section{Estimates for $s_0(D)$ and a characterization of the ball in ${\Bbb C}^n$}

\subsection{Characterization of the ball by determining
sets.}

This section is devoted to the proof of the following theorem.

\bT \label{ch}Let $D$ be a bounded domain in ${\Bbb C}^n$. Then $s_0(D)=n+1$ if and only if $D$ is biholomorphic to the unit ball $B^n$ in ${\Bbb C}^n$.\eT

To verify the estimate $s_0(B^n)=n+1$ we need to prove that no $n$ points in $B^n$ form a determining set for $Aut(B^n)$. (This was done in [FK2], we repeat it here for completeness). Consider $n$ arbitrary points 
$(p_0,p_1,...,p_{n-1})$, where $p_i\in B^n$ for $i=0,...,{n-1}$. Consider $g\in Aut(B^n)$ such that 
$g(p_0)=0$.  Consider now $n-1$ vectors $g(p_i)$, and the complex 
linear space $\pi $ spanned by these vectors.  Since $\dim (\pi )\leq n-1$, there is a rotation $f\in Aut(B^n)$ that is not the identity and 
keeps all the points of $\pi $ fixed. Now the automorphism 
$h=g^{-1}\circ f\circ g$ $\in \ Aut(B^n)$ is not the identity and it 
fixes all $n$ points $(p_0,p_1,...,p_{n-1})$. We proved that $W_n(B^n)=\emptyset$, so $s_0(B^n)=n+1$.
\newline
\newline The rest of this section will be devoted to the proof that $s_0(D)=n+1$ implies that $D$ is biholomorphic to the unit ball.
\newline

If $H$ is (isomorphic to) a subgroup of the unitary group $U(n)$, let $k(H)$ denote the least number $k$ of vectors $u_1,\dots, u_k$ such that if $h\in H$ and if $h(u_j)=u_j$ for $j=1,\dots,k$ then $h=id$. For $z\in D$ the isotropy group $Aut_z(D)$ is isomorphic to the group of its differentials at $z$, and these differentials are unitary with respect to the Bergman inner product on the tangent space $T_z(D)$. So $Aut_z(D)$ is isomorphic to a subgroup of $U(n)$.

\bL
For a bounded domain $D$ in ${\Bbb C}^n$, $s_0(D)\le 1+\min\{k(Aut_x(D)): x\in D\}$.
\eL

\begin{pf} Choose $z\in D$ so that $k(Aut_z(D))=\min\{k(Aut_x(D)): x\in D\}$. Denote that number by $k$. Let $u_1,\dots, u_k$ be vectors in $T_z D$ such that if $h\in Aut_z(D)$ and if $dh(z)(u_j)=u_j$ for $j=1,\dots,k$ then $dh=id$ (hence $h=id$). For each $u_j$, let $z_j$ be a point on the geodesic through $z$ in the direction $u_j$, so close to $z$ that the geodesic is the unique length minimizing geodesic from $z$ to $z_j$. Let $f$ be an automorphism of $D$ fixing $z, z_1, \dots, z_k$. Then $df(z)$ fixes $u_1, \dots, u_k$. It follows that $df(z)=id$ and $f=id$. Therefore, $s_0(D)\le  1+\min\{k(Aut_x(D)): x\in D\}$.
\end{pf}

\bL
If $H$ is a subgroup of $U(2)$ and if $H$ is not transitive on $S^3$, then $k(H)\le 1$.
\eL
\begin{pf}Let $S=S^3$ be the unit sphere in ${\Bbb C}^2$. It suffices to show that the set of fixed points in $S$ of nontrivial elements of $H$ (that is, each of these points is a fixed point of at least one nontrivial element of $H$) is not equal to $S$. For $g, h\in U(2)$ and $x\in S$, $x$ is a fixed point of $h$ iff $g^{-1}x$ is a fixed point of $g^{-1}hg$. So, without any loss of generality we can replace $H$ with a subgroup of $U(2)$ conjugate to $H$.

The Lie algebra $Q$ of $U(2)$ consists of skew Hermitian matrices, so it has as a basis the following elements:
$$a=\begin{pmatrix} i&0\\ 0&0\end{pmatrix},\;\;
b=\begin{pmatrix}0&0\\ 0&i\end{pmatrix},\;\; 
c=\begin{pmatrix}0&-1\\1&0\end{pmatrix},\;\;
d=\begin{pmatrix}0&i\\ i&0\end{pmatrix}.$$
Their Lie bracket relations are 
$$[a,b]=0, [a,c]=-d, [a,d]=c, [b,c]=d, [b,d]=-c, [c,d]=-2a+ 2b.$$

If $\dim H=4$, then $H=U(2)$, which contradicts the hypothesis that $H$ is not transitive on $S$.

Suppose that $\dim H=3$. One can verify that the only 3-dimensional Lie subalgebra of $Q$ is spanned by $\{c, d, a-b\}$. Hence, the identity component $H_0$ of $H$ is $SU(2)$, again contradicting the hypothesis that $H$ is not transitive on $S$.

Now suppose that $\dim H=2$. One can verify directly that $U(2)$ does not have a subgroup of dimension 2 and rank 1. Thus, $H$ has rank 2.  
Up to conjugation, the identity component $H_0$ of $H$ is $T^2=\{diag(e^{i\alpha}, e^{i\beta}): \alpha, \beta\in {\Bbb R}\}$. Each component of $H$ is $gT^2$. 
If $h\in gT^2$ has a nonzero fixed vector, $h$ must satisfy $det(h-id)=0$. It follows that if $\dim H\le 2$ then the set $U$ of the nontrivial elements of $H$ that have a fixed point on $S$ has dimension $\le 1$, and for each $g\in U$, the set of fixed points of $g$ on $S$ has dimension 1. Thus the set $P$ of fixed points of nontrivial elements of $H$ has dimension $\le 2$. It follows that $P\ne S$. Therefore, $k(H)\le 1$.

\end{pf}

\bL
If $H$ is a subgroup of $U(n)$ with $n\ge 2$ and if $H$ is not transitive on $S^{2n-1}$ then $k(H)\le n-1$.
\eL
\begin{pf}
The case where $n=2$ is the previous lemma. Suppose $n>2$ and $H\subset U(n)$ is not transitive on $S^{2n-1}$. Choose $x, y\in S^{2n-1}$ so that no element of $H$ maps $x$ to $y$. Choose $z\in S^{2n-1}$ orthogonal to both $x$ and $y$. Let $S_1=\{v\in S^{2n-1}: (v, z)=0\}$, where $(v,z)=\sum v_j\overline z_j$, and let $H_1=\{g\in H: g(z)=z\}$, and $U_1=\{g\in U(n): g(z)=z\}$. Now $U_1\cong U(n-1)$, and $H_1$ is a subgroup of $U_1$. By the induction hypothesis, since $H_1$ is not transitive on $S_1$, $H_1$ has a determining set of $n-2$ vectors $\{w_1, \dots, w_{n-2}\}$. It follows that $\{z, w_1,\dots, w_{n-2}\}$ is a determining set for $H$. Therefore, $k(H)\le n-1$.
\end{pf}
The proof of our main theorem follows from the following 
\bL If $D$ is a bounded domain in ${\Bbb C}^n$, and $D\not\cong B^n$, then $s_0(D)\le n$.
\eL
\begin{pf} If $n=1$ the statement is obviously true. Assume $n\ge 2$. Let $z\in D$. Since $D\not\cong B^n$, $Aut_z(D)$ is not transitive on the directions at $z$, by the main result of [GK]. By the Lemma~1.4, $k(Aut_z(D))\le n-1$. It follows that $s_0(D)\le 1+ k(Aut_z(D))\le n$.
\end{pf}
{\it Remark}. For endomorphisms $H(D,D)$ we still have for the ball $\widehat{s}_0(B^n)=n+1
$. However in the case of endomorphisms there are domains not biholomorphic to the ball but with the same maximum possible value of $\widehat{s}_0$. Here are two examples.
\newline For $n=1$, $\widehat{s}_0(D)=\widehat{s}_0(B^1)=2$ for any bounded domain $D\subset \Bbb C$.
\newline We will show now that for $n=2$, for
the unit polydisc $\Delta ^2$,  $\widehat{s}_0(\Delta ^2)=\widehat{s}_0(B^2)=3$.

Indeed, consider any two distinct points  $p_1,p_2\in \Delta ^2$. Since the $%
Aut(\Delta ^2)$ is transitive, we can find an automorphism $g$, such that $%
g(p_1)=0$. Let $L$ be the complex line through the origin and $g(p_2)$. In
terms of the coordinates this line can always be described in one of the
forms: $z_2=\lambda z_1$, or $z_1=\lambda z_2$ where $|\lambda |\leq 1$. One
can check that the map $P:(z_1,z_2)\rightarrow (z_1,\lambda z_1)$ in the
first case, or $P:(z_1,z_2)\rightarrow (\lambda z_2,z_2)$ in the second case
will produce a holomorphic retraction of the polydisc, fixing  $g(p_1),g(p_2)$. Now the
map $g^{-1}Pg$ is a holomorphic retraction of $\Delta ^2\,$fixing $p_1,p_2$. Therefore $%
\widehat{s}_0(\Delta ^2)>2$. Since this number is also $\leq 3$, we conclude 
$\widehat{s}_0(\Delta ^2)=3$.

\subsection{An estimate for $s_0(D)$.}
Let $G$ be a subgroup of $Aut(D)$. By $s_0(D,G)$ we denote the minimum
number of distinct points in $D$ such that if $g\in G$, and $g$ fixes all
these points, then $g=id$. So, $s_0(D)=s_0(D,Aut(D)).$

\bT\label{es} Let $D$ be a bounded domain in ${\Bbb C}^n$, let $G$ be a subgroup of $Aut(D)$, and let $q=\dim G$. If $q\ge 1$, then $s_0(D, G)\le q$. If $q=0$, then $s_0(D, G)\le 1$.
\eT
\begin{pf} First we consider the case where $q\le 1$. Let $e$ denote the identity element of $G$, and let $Q=G\backslash \{e\}$. For each $g\in Q$, the set $\{x\in D: g(x)=x\}$ is an analytic set of $D$ of dimension $\le 2n-2$. The set $W_1:=\{(g,x)\in Q\times D: g(x)=x\}$ is an analytic set of $Q\times D$ of dimension $\le (2n-2)+q\le 2n-1< dim D$. Let $W$ denote the set of fixed points of nontrivial elements of $G$. Since $W=\pi(W_1)$, where $\pi: Q\times D\to D$ is the projection, and since $dim W_1< dim D$, we see that $W\ne D$. Therefore, $s_0(D, G)\le 1$. 

Now we assume that $q\ge 2$. There must be an orbit $Q$ of $G$ of positive dimension. Let $x\in Q$, and let $H:=G_x$ be the subgroup of $G$ consisting of elements $g$ satisfying $g(x)=x$. Then $\dim H< \dim G$. By induction hypothesis, $s_0(D, H)\le \dim G-1$. Therefore, $s_0(D, G)\le 1+s_0(D, H)\le \dim G$.
\end{pf}

\bC\label{in}  Let $D$ be a bounded domain in ${\Bbb C}^n$. If 
$\dim (Aut(D))\ge 1$, then $s_0(D)\le \dim (Aut(D))$. If $\dim (Aut(D))=0$, then $s_0(D)\le 1$.
\eC

\section{On topological properties of determining sets. }
\subsection{Determining sets $W_s(D)$ are open and dense}

Our aim in this section is to prove the following theorem.

\bT\label{op} Let $D$ be a bounded domain in ${\Bbb C}^n$ and $s\ge 1$. Then $W_s(D)\subset D^s$ is open;
if in addition $W_s(D)\neq \emptyset $, then $W_s(D)$ is dense in $D^s$.
\eT

The assertion that $W_s(D)$ is open and dense in $D^s$ was proved for some domains and $s\ge n+1$ in [FK1]. Using analytic methods of [Ca1],[Ca2], J.P. Vigu\'{e} (see [Vi1],[Vi2]) proved that $W_s(D)$ is open for all bounded domains and all $s$, and that it is dense for $s\ge n+1$. By using the Bergman metric on a bounded domain we are able to use differential geometry methods and the Lie group properties (see [GKM], [Kl], [MZ], also [BD],[FMP], [Ma]) to prove the above general theorem.   

First we introduce some notation. If $G$ is a subgroup of $\Aut(D)$, 
$W_{s,G}(D)$ denotes the set of $s$-tuples $(x_1,\dots, x_s)$, where $x_j\in D$, such that each element $g\in G$ satisfying $g(x_j)=x_j$ for $j=1,\dots,s$ has to be the identity.

We need the following lemma (Theorem~2.4 in [Ma]).

\bL Let $\Omega $ be a bounded domain in ${\Bbb C}^n$ containing the closure of the unit ball, and $G$
a compact Lie subgroup of $Aut(\Omega )$. Suppose that each $G$-orbit lies
in a ball of radius $1/2$. Then $G=\{id\}$.\eL

\bL Suppose that $D$ is a bounded domain in ${\Bbb C}^n$ and $G$ is a subgroup of $Aut(D)$. Then $W_{1,G}(D)$ is open in $D$.
\eL
\begin{pf} We need to consider only the case $W_{1,G}(D)\neq \emptyset $ .
Suppose that $x\in W_{1,G}(D)$, and there is a sequence of points in $D$, $x_k\rightarrow x$ such that $x_k\notin W_{1,G}(D)$ $\forall k$. Let $U$ be a
neighborhood of $x$, $\overline{U}\subset D$, then there is a positive $r$
such that for large enough $k$ the ball with the center in $x_k$ and of
radius $r$, $b(x_k,r)$ in the \textit{Bergman metric }compactly lies in
\textit{\ }$U$. The assertion that $x_k\notin W_{1,G}(D)$ means that the
subgroup $G_{x_k}$ of $G\,$fixing $x_k$ is not the identity. This subgroup is a
compact Lie subgroup of $G$, and also acts on $b(x_k,r)$ (since Bergman
metric is an invariant metric). Applying (a properly adjusted form of) Lemma 2.2, one concludes that there exists an $\varepsilon >0$, such that for large
enough $k$ one has an automorphism $g_k\in G_{x_k}$ and a point $y_k\in
b(x_k,r)\subset U$ such that the Eucledian distance $|g_k(y_k)-y_k|>%
\varepsilon $. One can now find a subsequence $\{k_j\}$ such that (1) $y_{k_j}\rightarrow y\in \overline{U}$, and (2) $g_{k_j}\rightarrow g\in G$. We conclude now that $g(x)=x$ and that $|g(y)-y|\geq \varepsilon $. This means that $x\notin W_{1,G}(D)$ which is a
contradiction. Therefore $W_{1,G}(D)$ is open in $D$.
\end{pf}
Let $\rho(\cdot, \cdot)$ denote the Bergman distance. Let $b(z,r)$ denote the Bergman ball with center $z$ and radius $r$. Let $\overline b(z,r)$ be the closure of $b(z,r)$ in $D$.

\bL Suppose that $D$ is a bounded domain in ${\Bbb C}^n$ and $G$ is a subgroup of $Aut(D)$. If $W_{1,G}(D)\neq \emptyset $ then $W_{1,G}(D)$ is dense in $D$.
\eL
\begin{pf} In this proof, let $W=W_{1,G}(D)$. Suppose that $W$ is not dense in $D$. Then the closure $K$ of $W$ in $D$ is not equal to $D$. Let $p$ be a boundary point of $K$ in $D$. Choose $r>0$ such that the closure of $b(p, 4r)$ in $D$ is compact and such that each pair of points of $b(p, 4r)$ is connected by a unique length-minimizing geodesic segment (in the Bergman metric). There exist points $z, w$ such that $\rho(z, p)<r$, $\rho(w,p)<r$, $w\in W$, and $z\not\in K$. Note that the orbit of $w$, $G(w)\subset W$. Let $Q=G(w)\cap \overline b(p, 4r)$. Then $Q$ is compact and $Q \subset W$. Let $u$ be a point of $Q$ nearest to $z$. Then $u$ is also a point of $G(w)$ nearest to $z$, and $R:=\rho(z, u)\le \rho(z,w)<2r$. Choose a point $y$ on the unique length-minimizing geodesic segment from $z$ to $u$ such that $y\not\in K$ and $y\ne z$. For each point $x$ of $G(w)$, we see that 
$$\rho(z, y)+\rho(y, x)\ge \rho(z, x)\ge \rho(z, u),$$
and that the two equalities hold simultaneously only if $x=u$. Hence, $\rho(z, y)+\rho(y, x)> \rho(z, u)=R$ for each $x\in G(w)$, $x\ne u$. It follows that $\rho(y,x)> R-\rho(z,y)=\rho(y,u)$ for each $x\in G(w)$, $x\ne u$. Therefore, $u$ is the unique point of $G(w)$ nearest to $y$. Since $y\not\in K$, there is a nontrivial $g\in G$ such that $g(y)=y$. Now $\rho(y, u)=\rho(g(y), g(u))=\rho(y, g(u))$ forces $g(u)=u$. Since $u\in W$, the map $g$ must be the identity, contradicting the fact that $g$ is not trivial. Therefore, $W_{1,G}(D)$ is dense in $D$.
\end{pf}

{\it Proof} of Theorem \ref{op}.

We need to prove this theorem only for $W_s(D)\neq \emptyset $. For $g\in Aut(D)$ let $Q_s(g)$ denote the mapping 
\[
Q_s(g):D^s\rightarrow D^s,\;\;Q_s(g)(z_1,\dots ,z_s)=(g(z_1),\dots ,g(z_s)).
\]
Let $G=\{Q_s(g):g\in Aut(D)\}$. Then $G\subset Aut(D^s)$, and $%
W_{1,G}(D^s)=W_s(D)$. By the previous lemmas, $W_s(D)$ is open and dense in $%
D^s$.

\subsection{Determining sets $\widehat{W}_s(D)$ that are not open.}

In [Vi2] it was proved that $\widehat{W}_s(D)$ is open in $D^s$ for any bounded taut domain in ${\Bbb C}^n$. 
Our aim in this section is to present an example in ${\Bbb C}^2$ of a bounded domain
such that the determining set $\widehat{W}_2(D)$ is not open in $D^2$. 

First we construct the set $D\subset {\Bbb C}^2$.

Denote $B_2=\{z\in {\Bbb C}^2||z|<2\},B_1=\{z\in {\Bbb C}^2||z|<1\},b_1=%
\{z=(z_1,z_2)\in {\Bbb C}^2||z_1|^2+|z_2+1|^2<(1-10^{-4})^2\},$ $%
b_2=\{z=(z_1,z_2)\in {\Bbb C}^2||z_1-0.02|^2+|z_2+1|^2<1\}$. Let $\Omega
=(B_2\backslash \overline{B_1})\cup b_1\cup (b_2\cap B_2)$.

Consider now pairs of points $p_j,q_j\in \Omega $ , $%
p_j=(-1.5,2^{-j}),q_j=(1.5,2^{-j})$, and bydiscs $U_j=\{z=(z_1,z_2)\in {\Bbb C
}^2||z_1|<1.5,|z_2-2^{-j}|<2^{-2^j}\};j=1,2,...$. And, finally domain $%
D=\Omega \cup \bigcup\limits_{j=1}^\infty U_j$.

Note the following properties:

1. $D$ is a connected domain, $D \subset B_2$.

2. The entire complex disc $\Delta _j=\{z=(z_1,z_2)\in B_2|z_2=2^{-j}\}
\subset D$ for all $j$.

3. $\lim_{j\rightarrow \infty }p_j=p_0=(-1.5,0),\lim_{j\rightarrow \infty
}q_j=q_0=(1.5,0);$ and for the disc $\Delta _0=\{z=(z_1,z_2)\in B_2|z_2=0\} $%
, $\Delta _0\cap D=\{z\in {\Bbb C}^2|1<|z_1|<2,z_2=0)$.

4. $U_i\cap U_j=\emptyset $ for $i\neq j$ for large enough $i,j$.

5. Let $l_0=(0,-10^{-4})\in \overline{b_1}$ be the ``tip'' of this ball, and the
point of $\overline{b_1}$ closest to the origin. Let $d$ denote the
Kobayashi distance in $B_2$ from the origin to $l_0$, $\overline{k}=%
\overline{k}(0,d)$ the closed Kobayashi ball with the center at the origin
and radius $d$. Then $\overline{k}\cap \overline{(b_1\cup b_2)}=\{l_0\}$.

We are going to show that for any $j$ the pair $(p_j,q_j)\notin \widehat{W}%
_2(D)$, but their limit $(p_0,q_0)\in \widehat{W}_2(D)$, which will prove
that $\widehat{W}_2(D)$ is not open in $D^2$.

{\it Statement} 1. For any $j$, there is a holomorphic retraction of $B_2$ (and therefore of $%
D\subset B_2$) onto $\Delta _j$, and since $(p_j,q_j)\in \Delta _j$, the
pair $(p_j,q_j)\notin \widehat{W}_2(D).$

To prove this one needs to use first an automorphism $g$ of $B_2$ to move $%
\Delta _j$ to $\Delta _0$, use the natural projection $P$ of the ball $B_2$
onto $\Delta _0$, and set the needed holomorphic retraction as $g^{-1}\circ P\circ g$.

{\it Statement} 2. Any holomorphic map $f:D\rightarrow D$ extends to a holomorphic map $F:B_2\rightarrow B_2.$

{\it Statement} 3. Any holomorphic map $F:B_2\rightarrow B_2$ that fixes our two
points $p_0,q_0$ will fix all the points of $\Delta _0$.

For proof see ([Vi2], ex. 1 in sec 4).
\newline
\newline Now let $f:D\rightarrow D$ be a holomorphic map that fixes our two points $%
p_0,q_0$. Its extension $F:B_2\rightarrow B_2$ will be an identity on $%
\Delta _0$.

{\it Statement} 4. $F(l_0)=l_0$.

Indeed, consider $K\,$the Kobayashi ball in $B_2$ with center at the origin that
coincides with the standard unit ball in ${\Bbb C}^2$. $K\cap D$ consists or
nonintersecting connected pieces, only one of which, namely $G=K\cap
(b_1\cup b_2)$ has a point $(0.02,0)$ on $\Delta _0$ as a limit point. Since
this point is fixed by $F$, and $F$ cannot increase the Kobayashi distance,
we conclude that $F(G)\subseteq G$.

Since it is also true that $F($ $\overline{k})\subseteq $ $\overline{k}$,
the only possible image for $l_0$ under $F$ is the point itself (see
property 5 above), so $F($ $l_0)=l_0$.

We now conclude our observation by pointing out that $F:B_2\rightarrow B_2$
has three fixed points $p_0,q_0,l_0$. By ([Vi2], ex. 1 in sec 4) $F$ is
the identity, so $f:D\rightarrow D$ is the identity, and therefore $%
(p_0,q_0)\in \widehat{W}_2(D)$.

\subsection{Determining sets $\widehat{W}_s(D)$ that are not dense.}

Our goal here is to present an example of a domain $D$ such that
for any $s$ the determining set $\widehat{W}_s(D)$ is not dense in $D^s$.

J.-P. Vigu\'{e} (see [Vi2], ex. 2 in sec 4) has provided such an example for $
D=\Delta ^2=\{z=(z_1,z_2)\in {\Bbb C}^2||z_1|<1,|z_2|<1\}$ the polydisc in $
{\Bbb C}^2$, and $s=3$. For completeness, using the same idea, we provide here an example for $\Delta ^2$ and any $s\geq 3$.

Fix $s$ points $(A_j,0)$, $A_j=2^{-j},j=1,...,s$. The set $T$ of these
points is a point in $(\Delta ^2)^s$. $T\notin \widehat{W}_s(\Delta
^2)$, since $(z_1,z_2)\mapsto (z_1,0)$ is a holomorphic retraction of $\Delta ^2$,
fixing all these points. Let $\widetilde{T}=\{(a_j,b_j),j=1,...,s\}\in
(\Delta ^2)^s$ be any $\delta >0$ perturbation of $T$. So, $\sum\limits_{j=1}^{s}(|A_j-a_j|^2+|b_j|^2)<\delta ^2$. We will show
that if $\delta $ is small enough $\widetilde{T}\notin \widehat{W}_s(\Delta
^2)$.

Consider the Lagrange interpolation polynomial $\varphi (w)=\sum\limits_{j=1}^sb_j\prod\limits_{i\neq
j}\frac{(w-a_i)}{(a_j-a_i)}$. One can verify that if $\delta >0$ is small
enough (say $\delta <4^{-s^2}$) then $|\varphi (w)|<1$ if $|w|<1$. Now the
map $(z_1,z_2)\mapsto (z_1,\varphi (z_1))\,$is a holomorphic retraction of the unit
polydisc that has $\widetilde{T}$ in its set of fixed points. Therefore $%
\widetilde{T}\notin \widehat{W}_s(\Delta ^2)$.

As a remark we note that using this idea one can construct many such
examples. Moreover the following theorem holds.

\bT Consider the topological space $\check{D}_n$ of all bounded domains in ${\Bbb C}^n$ with the topology induced by the Hausdorff distance between boundaries
of the domains. Let $\breve{G}_n\subset \check{D}_n$ be such that if $D\in 
\breve{G}_n$ then $\widehat{W}_s(D)$ is not dense in $D^s$ for all $s\geq 1$%
. Then $\breve{G}_n$ is dense in the topological space $\check{D}_n$.
\eT

\end{document}